\documentclass{proc-l}

\usepackage{psfig}

\newtheorem{theorem}{Theorem}[section]
\newtheorem{lemma}[theorem]{Lemma}

\theoremstyle{definition}

\newtheorem{corollary}[theorem]{Corollary}
\newtheorem{proposition}[theorem]{Proposition}

\newcommand{\cS}{\mathcal{S}_{2,\infty}}
\newcommand{\cL}{\mathcal{L}_{fr}}

\newcommand{\bZ}{{\mathbb Z}}
\newcommand{\bC}{{\mathbb C}}
\newcommand{\bQ}{{\mathbb Q}\cup\{\frac{1}{0}\}}
\newcommand{\bN}{{\mathbb N}}
\newcommand{\bH}{{\mathbb H}^2}
\newcommand{\bR}{{\mathbb R}}
\newcommand{\g}{SL_2(\bC)}

\newcommand{\lcr}{\raisebox{-5pt}{\mbox{}\hspace{1pt}
                  \psfig{file=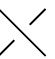}\hspace{1pt}\mbox{}}}
\newcommand{\ift}{\raisebox{-5pt}{\mbox{}\hspace{1pt}
                  \psfig{file=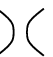}\hspace{1pt}\mbox{}}}
\newcommand{\zer}{\raisebox{-5pt}{\mbox{}\hspace{1pt}
                  \psfig{file=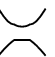}\hspace{1pt}\mbox{}}}
\newcommand{\smst}[1]{\makebox[0pt]{\scriptsize{$#1$}}}

\theoremstyle{remark}

\numberwithin{equation}{section}

\begin{document}

\title[Kauffman bracket skein quantizations]{multiplicative structure
of kauffman bracket skein module quantizations}

% author one information
\author[d. bullock]{Doug Bullock} 
\address{Department of Mathematics,
University of Maryland, College Park, MD 20742, USA}
\curraddr{} 
\email{bullock@math.gwu.edu} 
\thanks{Supported by an NSF postdoctoral fellowship.}

% author two information
\author[j. h. przytycki]{J{\'o}zef H. Przytycki}
\address{Department of Mathematics,
The George Washington University, Washington, DC 20052, USA}
\curraddr{}
\email{przytyck@math.gwu.edu}
\thanks{Supported by USAF grant 1-443964-22502 while visiting the  
 Mathematics Department, U.C. Berkeley.}

\date{}

% at present the "communicated by" line appears only in ERA and PROC
\commby{}

\dedicatory{}

\begin{abstract}
We describe, for a few small examples, the Kauffman bracket skein
algebra of a surface crossed with an interval.  If the surface is a
punctured torus the result is a quantization of the symmetric algebra
in three variables (and an algebra closely related to a cyclic
quantization of $U(\mathfrak{so}_3$).  For a torus without boundary we
obtain a quantization of ``the symmetric homologies" of a torus
(equivalently, the coordinate ring of the $SL_2(\bC)$-character
variety of $\bZ\oplus\bZ$).  Presentations are also given for the four
punctured sphere and twice punctured torus.  We conclude with an
investigation of central elements and zero divisors.
\end{abstract}

\maketitle

\section{Introduction}

The {\em Kauffman bracket skein module}, $\cS(M)$, is defined as
follows \cite{H-P-1,Pr-1}: Assume that $M$ is an oriented 3-manifold,
that $\cL$ denotes unoriented framed links in $M$ (including the empty
link), and that $R$ is any commutative ring with unit.  With $A$
invertible in $R$, let $S_{2,\infty}$ be the submodule of $R\cL$
generated by all expressions of the form
\[\lcr-A\zer-A^{-1}\ift \quad \text{and} \quad \bigcirc+A^2+A^{-2}.\]
The diagrams in each relation indicate framed links that can be
isotoped to identical embeddings except within the neighborhood shown,
where framing is vertical.  Set $\cS(M;R,A)= R\cL/S_{2,\infty}$.
Notation is shortened for two special cases: $\cS(M) =
\cS(M;\bZ[A^{\pm 1}],A)$ and $\mathcal{S}(M) = \cS(M;\bZ,-1)$.

The following standard results about skein modules will be useful
later on.  Only for (1) and (6) are the proofs more than elementary
exercises.  These two may be found in \cite{skeinalg} and \cite{Pr-1}.

\begin{proposition}\mbox{}
\begin{enumerate}
\item 
(Universal coefficients property).  Let $r: R \to R'$ be a
homomorphism of rings (commutative with 1), making $R'$ into an $R$
module. The identity map on $\cL$ induces an isomorphism of $R'$
(and $R$) modules:
$$ \cS(M;R',r(A)) \cong \cS(M;R,A)\otimes_R R' .$$ 

\item 
An embedding of 3-manifolds $f:M \to N$ induces the homomorphism of
skein modules $f_*:\cS(M;R,A) \to \cS(N;R,A)$. This leads to a functor
from the category of 3-manifolds and embeddings to the category of $R$
modules.  If $N$ is obtained from $M$ by adding 2- and 3-handles, then
$f_*$ is an epimorphism.

\item 
If $M=F\times I$ for an oriented surface $F$ and an interval $I$, then
$\cS(M;R,A)$ is an $R$ algebra with $\emptyset$ as a unit element and
$L_1\cdot L_2$ defined by placing $L_1$ above $L_2$.  This
multiplication depends on the product structure of $M$, so we use the
notation $\cS(F;R,A)$.  An embedding of oriented surfaces $f:F \to F'$
induces the homomorphism of skein algebras $f_*:\cS(F;R,A) \to
\cS(F';R,A)$. This leads to a functor from the category of surfaces and
embeddings to the category of $R$ algebras.

\item 
If $A=-1$ then, for any $M$, $\cS(M,R,-1)$ is an $R$ algebra:
$L_1\cdot L_2$ is defined to be the disjoint union of links. The
algebra depends only on $\pi_1(M)$.  In particular, if $f:M \to N$ is
a homotopy equivalence then $f_*: \cS(M,R,-1) \to \cS(N,R,-1)$ is an
isomorphism of algebras.

\item 
The skein module $\cS(F;R,A)$ has a basis consisting of links on $F$
without contractible components (but including an empty link).

\item 
The skein algebra $\cS(\raisebox{-2pt}{\psfig{figure=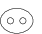}};R,A)$
is isomorphic to $R[x,y,z]$ where the indeterminates are the boundary
curves.
\end{enumerate}
\end{proposition}

Let $F_{g,n}$ denote an oriented surface of genus $g$ with $n$
boundary components.  We will compute $\cS(F_{g,n};R,A)$ for $(g,n)
\in \{(1,0),(1,1),(1,2),(0,4)\}$.  One can intuit a presentation from
existing computations of the $SL_2(\bC)$-character ring of a free
group \cite{Hor} and its relationship with $\cS(F_{g,n};\bC,-1)$
\cite{B,PS}.  Verification is then a matter of confirming relations
and constructing a basis.  Often we will check linear independence by
specializing to a simpler module:

\begin{lemma} 
Let $(u)$ be a principal ideal in a Noetherian ring $R$, and let $\pi:
R \to R/(u)$ be the natural epimorphism.  Suppose $\{v_i\}$ is a
subset of a torsion free $R$ module $S$.  If $\{v_i \otimes 1\}$ are
linearly independent in the $R/(u)$ module $S \otimes_R R/(u)$, then
$\{v_i\}$ are linearly independent in $S$.
\end{lemma}
 
\begin{proof}
Suppose $\sum a_iv_i =0$ with some $a_j \neq 0$.  Since $R$ is
Noetherian, there is a maximal $k$ such that $u^k$ divides each $a_i$.
Let $a_i = u^k a'_i$.  No torsion in $S$ implies $\sum a'_iv_i =0$, so
$\sum (a'_i+(u))(v_i \otimes 1) = 0$ in $S \otimes_R R/(u)$.  Thus
each $a'_i \in (u)$, contradicting maximality of $k$.
\end{proof}

In particular, one sees that a linearly independent set of links in
$\mathcal{S}(M)$ is also linearly independent in $\cS(M)$.

\section{Skein algebras of $F_{1,1}$ and $F_{1,0}$}

For the rest of this paper let us agree on the shorthand nomenclature
{\em curve} for a simple closed curve (up to isotopy) in a surface $F$,
and for the knot in $F \times I$ created by framing it vertically.
Suppose that $x_1$ and $x_2$ are curves in $F_{1,1}$ that intersect
once.  Applying a skein relation to the single crossing of the link
$x_1x_2$ resolves it into
\begin{equation}
x_1x_2 = A x_3 + A^{-1} z,
\end{equation}
where $z$ and $x_3$ are the unique curves that meet both $x_1$ and
$x_2$ once.  For a unit $u \in R$ and $a$ and $b$ in an $R$ algebra,
let $[a,b]_u$ denote the deformed commutator $uab -u^{-1}ba$.  Also,
let $\delta = A^2 - A^{-2}$.  We use $R\langle\{g\}\;|\;\{r\}\rangle$
to denote the free $R$ algebra in non-commuting variables $\{g\}$
modulo the ideal generated by $\{r\}$.  (Often the relations $r$ will
be written as equations.)

\begin{theorem}
With curves $x_i$ as above, $\cS(F_{1,1};R,A)$ is presented as
\[R\langle x_1,x_2,x_3\ |\ [x_i,x_{i+1}]_A  = \delta x_{i+2} \rangle,\]
where $i=1,2,3$ and subscripts are interpreted modulo three.  
\end{theorem}

\begin{proof}
The skein algebra is generated by curves in the surface (Proposition
1(5)), which, except for the boundary, may be identified with slopes
in $\bQ$.  These in turn can be organized as the vertices in a
tessellation of the upper half space model of $\bH$ by ideal
triangles---their sides are geodesics connecting slopes of curves that
cross once.  The product of two curves that meet once resolves, with
invertible coefficients, into two other curves that meet each of the
original pair in a single point, e.g.\ Equation (1).  Hence, one can
express the vertices of any triangle in terms of an adjacent one.  It
follows that every non-boundary parallel curve is generated by
$\{x_i\}$.  The boundary, $\partial$, is obtained from the resolution
\begin{equation}
z x_3 = A^2 x_1^2 + A^{-2} x_2^2 - A^2 - A^{-2} + \partial,
\end{equation}
where $z$ is as is Equation (1).

The commutation relations follow easily by resolving the links
$x_ix_j$. Any other relation among $x_1$, $x_2$ and $x_3$ reduces,
using the commutators, to an expression in $\{x_1^ix_2^jx_3^k\}$.  It
suffices, then, to show that this set is a basis.  Since
$\cS(F_{1,1};R,A)$ is free, the universal coefficients property allows
us to restrict to $\cS(F_{1,1})$.  By Lemma 1, linear independence may
be checked in $\mathcal{S}(F_{1,1})$, where it follows from
Proposition 1(6).
\end{proof}

The presentation in Theorem 1 is best described as a {\em cyclic
deformation} of $R[x_1,x_2,x_3]$.  There is a cyclic deformation of
$U(\mathfrak{so}_3)$ (Zachos terms it the cyclically symmetric Fairle
rotation \cite{Od,Za-1}) given by $U_A(\mathfrak{so}_3) = R\langle
y_1,y_2,y_3\ |\ [y_i,y_{i+1}]_A = y_{i+2} \rangle$.  This is related
to the skein algebra as follows.

\begin{corollary}
There is a well defined map from $\cS(F_{1,1};R,A)$ to
$U_A(\mathfrak{so}_3)$ sending $x_i$ to $\delta y_i$.  It is injective
if and only if $\delta$ is not a zero divisor and surjective if and
only if $\delta$ is invertible.
\end{corollary}

This is particularly interesting if $R = \bC$, for it says there are
two families of algebras parameterized by $A$---indistinguishable at
generic $A$ but having different fibers at $A = \pm 1$, where one family
has $U({\mathfrak so}_3)$ and the other has $\bC [x_1,x_2,x_3]$.

It has also been observed that, generically, the algebra
$\cS(F_{1,1}\times I;R,A)$ cannot be generated by fewer than three
elements.  It is clear, however, that one of the generators may be
eliminated if $\delta$ is invertible.

Finally, we obtain a presentation of the skein algebra of a torus by
adjoining one relation to the presentation in Theorem 1.

\begin{theorem}
$\cS(F_{1,0};R,A)$ is the quotient of $\cS(F_{1,1};R,A)$ by the
principle ideal $(A^2 x_1^2 + A^{-2} x_2^2 + A^2 x_3^2 - A x_1x_2x_3
- 2A^2 - 2A^{-2}).$
\end{theorem}

\begin{proof}
Embedding $F_{1,1}$ into $F_{1,0}$ maps $\cS(F_{1,1};R,A)$ onto
$\cS(F_{1,0};R,A)$, forcing $\partial = - (A^2 + A^{-2})$.  Any
relation that is not a multiple of $\partial +A^2 +A^{-2}$ would imply
a non-trivial relation among the standard basis elements in
$\cS(F_{1,0};R,A)$.  Therefore, we need only express $\partial$ in the
generators $x_1,x_2,x_3$. This is accomplished by eliminating $z$ from
Equations (1) and (2).
\end{proof}

\section{The Skein Algebra of $F_{0,4}$}

Paralleling Section 2, chose distinct, non-boundary parallel curves
$x_1$ and $x_2$ on $F_{0,4}$ that intersect in two points.  Then
compute
\begin{equation}
x_1 x_2 = A^2 x_3 + A^{-2} z + \text{boundary curves},
\end{equation}
where $x_3$ and $z$ are the unique curves whose minimal intersection
with each of $x_1$ and $x_2$ is a pair of points.  The similarity with
Equation (1) should call to mind the triangular tessellation of $\bH$
from the previous section.  The vertices can be matched to the
non-boundary parallel curves on $F_{0,4}$ so that each side of a
triangle joins curves that meet in two points.  Moreover, Equation (3)
holds with the variables replaced by the vertices of any adjacent pair
of triangles, provided the product term is the pair of shared
vertices.

We can extend $R$ to $\bar R$ inside $\cS(F_{0,4};R,A)$ by adjoining
the boundary components of $F_{0,4}$.  Since these are central,
$\cS(F_{0,4};R,A)$ is an algebra over the commutative ring $\bar R$.
Links in the surface with no trivial or boundary parallel components
form and $\bar R$-basis.

Suppose that $x_1$ separates the boundary curves $a_1$ and $a_2$ from
the boundary curves $a_3$ and $a_4$.  Set $p_1 = a_1a_2 + a_3a_4$.
Define $p_2$ and $p_3$ similarly.  Let $q = a_1 a_2 a_3 a_4 + a_1^2 +
a_2^2 + a_3^2 + a_4^2$.  One may check, by direct resolution or by
eliminating $z$ from a resolution of $zx_3$, that
\begin{align}
A^2x_1x_2x_3 &= A^4x_1^2 + A^{-4}x_2^2 + A^4x_3^2 + A^2p_1x_1 +
A^{-2}p_2x_2 + A^2p_3x_3 \\ &+ q - (A^2+A^{-2})^2. \nonumber
\end{align}

The similarity with the relation in $\cS(F_{1,0};R,A)$, along with the
$\bH$ parameterization of curves in $F_{1,0}$ and $F_{0,4}$, suggests
that one can pass between their skein algebras by sending $A$ to $A^2$
while requiring the boundary curves to satisfy\footnote{This was
originally motivated by the double branched cover of the torus over
the sphere, with four branch points.}
\begin{equation}
p_i = 0 \quad \text{and} \quad q + \delta^2 = 0.
\end{equation}
In any $R$ there is at least one solution given by $a_1 = a_2 = a_3 =
-a_4 = A+A^{-1}$.  Let $J$ be an ideal in $\cS(F_{0,4};R,A)$ generated
by a solution of Equation (5).

\begin{theorem}
With notation as above, $\cS(F_{0,4};R,A)$ is presented as
\[\bar{R} \langle x_1,x_2,x_3\ |\ [x_i,x_{i+1}]_{A^2} = (A^4 - A^{-4})
x_{i+2} - \delta p_{i+2}, {\rm Equation (4)}\rangle.\]
\end{theorem}

\begin{proof}
The curves $x_i$ generate as in the proof of Theorem 1.  Commutators
and Equation (4) are a matter of direct computation.  The relations
imply that the set $\{x_1^ix_2^jx_3^k\;|\; ijk=0\}$ spans the module,
so we need only show they are linearly independent.  As before, it
suffices to consider the case $(R,A)=(\bZ,-1)$.  Here, formal
identification of variables gives an isomorphism $\mathcal{S}(F_{1,0})
\cong \mathcal{S}(F_{0,4})/J$.  One may check that
$\{x_1^ix_2^jx_3^k\;|\;ijk=0\}$ is a basis for $\mathcal{S}(F_{1,0})$.
In particular, the same set of monomials is linearly independent in
$\mathcal{S}(F_{0,4})/J$.  Finally, Lemma 1 (used four times) proves
they are linearly independent in $\mathcal{S}(F_{0,4})$.
\end{proof}

\begin{corollary}
One obtains a presentation of $\cS(F_{0,4};R,A)$ over $R$ from the
one in Theorem 3 by adding generators $a_1, \ldots, a_4$, and by
adding commutation relations to make them central.
\end{corollary}

\begin{proof}
The set $\{x_1^ix_2^jx_3^ka_1^pa_2^qa_3^ra_4^s\;|\;pqrs=0\}$ is seen
to be a basis by working over $(R,A) = (\bZ,-1)$. 
\end{proof}

\section{The skein algebra of $F_{1,2}$}

This algebra is harder to describe because there are far too many
simple closed curves on the surface.  However, since $F_{1,2}$ and
$F_{0,4}$ are homotopy equivalent, we are able to guess and verify a
presentation of $\cS(F_{1,2};R,A)$.

\begin{figure}
\mbox{}\hfill
\psfig{file=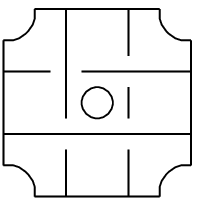}\hspace{-48pt}\raisebox{21pt}{\smst{x_1}}
\hspace{15pt}\raisebox{7pt}{\smst{y_1}}\hspace{7pt}\raisebox{48pt}{\smst{y_2}}
\hspace{12pt}\raisebox{33pt}{\smst{x_2}}\hfill\mbox{}
\caption{}
\label{curves}
\end{figure}

Choose curves $x_1$, $x_2$, $y_1$ and $y_2$ on $F_{1,2}$ as indicated
in Figure \ref{curves}.  (Edges are identified top-to-bottom and
left-to-right.)  Let $a$ denote the boundary curve in the center of
the figure and define $\bar R = R[a] \subset \cS(F_{1,2};R,A)$.
Define curves $z_1$ and $z_2$ by the resolutions $x_1y_1 = Az_1 +
A^{-1}w_1$, and $y_1x_2 = Az_2 +A^{-1}w_2$.  

There is a homotopy equivalence between $F_{1,2}$ and $F_{0,4}$ given
by retracting either surface to the spine in Figure \ref{spines} and
then identifying it with the spine in the other surface.  The induced
map turns Equation (4) into a relation in $\mathcal{S}(F_{1,2})$.  By
experimenting with coefficients we discovered the following relation
in $\cS(F_{0,4};R,A)$.
\begin{align}
A^2az_2z_1 &= A^2a^2+A^{-2}z_2^2+A^6z_1^2 +(y_1y_2+A^4x_1x_2)a \\
&- (A^{-1}x_1y_2 + A^{-1}x_2y_1)z_2 - (Ax_2y_2+A^5x_1y_1)z_1
\nonumber\\
&+x_2y_1x_1y_2+A^6x_1^2+A^2x_2^2+A^2y_1^2+A^{-2}y_2^2-A^2(A^2+A^{-2})^2
\nonumber
\end{align}
Let $C(t_1,t_2,t_3)$ denote the set of cyclic commutators,
$\{[t_i,t_{i+1}]_A=\delta t_{i+2}\;|\;i=1,2,3\}$, where subscripts are
taken modulo three.

\begin{figure}
\mbox{}\hfill\psfig{file=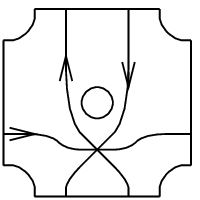}\hspace{-43pt}\raisebox{43pt}
{\smst{\alpha}}\hspace{18pt}\raisebox{43pt}{\smst{\beta}}\hspace{12pt}
\raisebox{23pt}{\smst{\gamma}}\hfill
\psfig{file=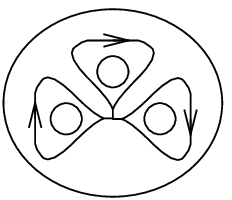}\hspace{-59pt}\raisebox{30pt}{\smst{\gamma}}
\hspace{23pt}\raisebox{49pt}{\smst{\alpha}}\hspace{25pt} 
\raisebox{30pt}{\smst{\beta}}\hfill\mbox{}
\caption{}
\label{spines}
\end{figure}

\begin{theorem}
$\cS(F_{2,1};R,A)$ is presented\,\footnote{A slightly different
presentation was quoted in the survey article \cite{Pr-2}} as
\begin{align*}
\bar{R} \langle x_i, y_i, z_i\;|\; 
& [x_1,x_2]=0, [y_1,y_2]=0, [z_1,z_2]=\delta (x_2y_2 -x_1y_1), \\
& C(x_1,y_1,z_1), C(x_2,y_2,z_1), C(z_2, y_1,x_2), C(z_2, y_2, x_1), 
 {\rm Equation (6)} \rangle.
\end{align*}
\end{theorem}

\begin{proof}
We present over $\bar R$ only to save the trouble of writing out
trivial commutators; in the proof we work $R$-linearly.  Given a
handle decomposition of $F_{1,2}$, there is a generating set for
$\cS(F_{1,2};R,A)$ given in \cite{Bu}.  It is possible to choose the
decomposition so that this set is $\{x_1,x_2,y_1,y_2,z_1,w_2,a\}$.
Since $w_2 = Ay_1x_2-A^2z_2$ we can replace it with $z_2$.  

Relations are checked by brute force.  Equation (6) and commutators
indicate that $\{a^iz_2^jz_1^kx_2^py_1^qy_2^rx_1^s\;|\;pqrs=0\}$ is a
spanning set.  Thus we need only show it is a basis in the case
$(R,A)=(\bZ,-1)$.  This follows from the homotopy equivalence in
Figure \ref{spines} and the proof of Corollary 3.
\end{proof}

\section{Central elements and zero divisors}

One can prove that $\cS(F)$ has no zero divisors and that its center
is the subalgebra generated by boundary curves: All that is needed is
a system for describing the standard basis so that the product of two
basic elements is conveniently expressed in this language. In
practice, one assigns a complexity to basis elements so that lower
order terms can be ignored.  In this section we will explicitly lay
out these ingredients for $F_{1,0}$, $F_{1,1}$ and $F_{0,4}$.

\subsection{Toral examples}

Links in $F_{1,0}$ (with no trivial components) are denoted by points
in $\bZ \times \bZ$, taking a point and its negative to represent the
same link.  The link $v=(p,q)$ is $\text{gcd}(p,q)$ copies of the
slope $p/q$ curve.  The complexity of $(p,q)$ is $p^2+q^2$.  The same
conventions apply to $F_{1,1}$, treating the boundary curve as a
scalar.  The notation $|v \cap w|$ means the minimal geometric
intersection number of $v$ and $w$ as links.  

\begin{lemma}
Choose $v$ and $w$ in $\bZ \times \bZ$.  The products $vw$ and $wv$
are expressed in the standard basis using only links in the
parallelogram spanned by $\{\pm(v\pm w)\}$.  The two distinct corners
occur with coefficients $A^{\pm|v \cap w|}$.
\end{lemma}

\begin{proof}
For $v = (n,0)$, induct on $n$.  The arbitrary case transforms to this
one using an element of $PSL_2(\bZ)$.
\end{proof}

\begin{theorem}
If $A$ is not a root of unity and $R$ has no zero divisors, then
the center of $\cS(F_{1,0};R,A))$ is $R$ and the center of
$\cS(F_{1,1};R,A))$ is $R[\partial]$.
\end{theorem}

\begin{proof}
Let $\alpha = \sum r_iv_i$ with each $r_i \neq 0$ in $R$ and each
$v_i$ a distinct, non-empty link in $F_{1,0}$.  Choose $w$ that is not
parallel to any $v_i$.  Suppose $u$ is a maximal complexity link
occurring in the resolution of $w\alpha - \alpha w$ into the standard
basis.  The resolution of each $r_i(wv_i-v_iw)$ lies in a distinct
parallelogram, so $u$ must be a vertex that is contained in exactly
one of them.  It cannot, therefore,  cancel with any other terms of the
resolution.  Since its coefficient is $\pm r_i (A^{|v_i \cap w|} -
A^{-|v_i \cap w|})$, $\alpha$ is not central.  The same proof works for
$F_{1,1}$, treating $\partial$ as a scalar.
\end{proof}

For zero divisors we will not be able to separate the resolution into
distinct parallelograms.  Instead we use an elementary fact about
vectors in $\bR^n$.

\begin{lemma}
Suppose $v_1$, $v_2$, $w_1$ and $w_2$ are vectors in $\bR^n$ with
$v_1+w_1 =v_2+w_2$.  If $v_1 \neq v_2$ then one of $v_1 + w_2$ or $v_2
+ w_1$ is longer than $v_1+w_1$.
\end{lemma}

\begin{proof}
Rotate so that $v_1+w_1$ lies in the first coordinate and assume that
each $||v_i + w_j|| \leq ||v_1 + w_1||$.  Use both conditions to force
the first coordinates of $v_1$ and $v_2$ agree, as well as those of
$w_1$ and $w_2$.  Agreement in the other coordinates follows easily.
\end{proof}

\begin{theorem}
$\cS(F_{1,0};R,A)$ and $\cS(F_{1,1};R,A)$ have zero divisors only if $R$
does.
\end{theorem}

\begin{proof}
Choose $\alpha = \sum r_i v_i \in \cS(F_{1,0};R,A)$ with each $v_i$
distinct and each $r_i \neq 0$.  Similarly, choose $\beta = \sum s_j
w_j$.  Suppose that $u$ has maximal complexity in the resolution of
$\alpha \beta = \sum_{i,j} r_is_j (v_iw_j)$.  Reordering if necessary,
we may assume that $u$ is one of $v_1 \pm w_1$.  If $u$ appears in the
resolution of any $v_iw_j$ other than $v_1w_1$, then Lemma 3 forces a
term of either $v_1w_j$ or $v_iw_1$ to have greater complexity.  Thus
$\alpha \beta$ expressed in the standard basis contains the non-zero
term $r_1s_1A^{\pm |v_1 \cap w_1|}u$.  The proof for $F_{1,1}$ is
similar, noting that $R[\partial]$ has no zero divisors.
\end{proof}

\subsection{The planar example}

Decorate $F_{0,4}$ with arcs $\alpha$, $\beta$ and $\gamma$ as in
Figure \ref{arcs}.  A link in $F_{0,4}$ with no trivial components and
no component parallel to an inner boundary curve defines a point
$(a,b,c) \in \bN\,^3$ as follows.  Isotope it to meet $\alpha \cup
\beta \cup \gamma$ minimally.  The number of arcs joining $\alpha$ to
$\beta$ is $a$; the number from $\beta$ to $\gamma$ is $b$, and from
$\gamma$ to $\alpha$, $c$.  A point corresponds to a link if its
coordinates are congruent modulo two.  There are two useful
complexities for $(a,b,c)$: $a + b + c$ and $a^2 + b^2 + c^2$.
Throughout this section we assume the three inner boundary components
are scalars in $\cS(F_{0,4};R,A)$.  Note that the extended scalars
have zero divisors only if $R$ does.

\begin{figure}
\mbox{}\hfill
\psfig{file=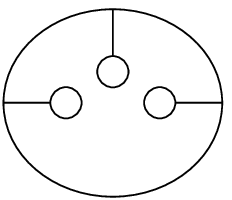}\hspace{-56pt}\raisebox{31pt}{\smst{\beta}}
\hspace{25pt}\raisebox{47pt}{\smst{\alpha}}\hspace{17pt} 
\raisebox{31pt}{\smst{\gamma}}\hfill\mbox{}
\caption{}
\label{arcs}
\end{figure}

\begin{lemma}
Using either complexity, \[(a,b,c) (a',b',c') =
A^{(ab'-ba'+bc'-cb'+ca'-ac')/2} (a+a',b+b',c+c') +\mathcal{O},\] where
$\mathcal{O}$ denotes lower complexity terms.
\end{lemma}

\begin{proof}
Assume for now that $a > c$.  Position the link so that near $\alpha$
it looks like Figure \ref{resolution}, in which a number next to a
strand indicates so many parallel copies.  If every crossing in Figure
\ref{resolution} is smoothed in the vertical direction the resulting
coefficient is $A^{(a-c)(a'+c')/2}$.  If $c < a$ the link should be
isotoped to look like a left-to-right reflection of Figure
\ref{resolution}, but vertical resolution of every crossing yields the
same coefficient.  In either case, a resolution including even one
horizontal smoothing allows an isotopy that removes two points of
intersection with $\alpha$ (possibly by creating a scalar boundary
component.)

\begin{figure}
\mbox{}\hfill
\psfig{file=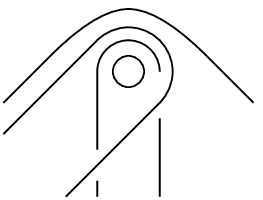}
\hspace{-76pt}\smst{\frac{1}{2}(a-c)}
\hspace{-3pt}\raisebox{33pt}{\smst{c}}
\hspace{56pt}\smst{a'+c'}\hfill\mbox{}
\caption{}
\label{resolution}
\end{figure}

Similarly position $(a,b,c)(a',b',c')$ near $\beta$ and $\gamma$.  The
rest of the link can be drawn with no more crossings.  Therefore the
complete resolution contains a term $A^{(ab'-ba'+bc'-cb'+ca'-ac')/2}
(a+a',b+b',c+c')$. Moreover, any other link in the resolution will be
some $(x,y,z)$ with $x \leq a$, $y \leq b$, $z \leq c$, and at least
one inequality strict.
\end{proof}

\begin{theorem}
If $R$ has no zero divisors and $A$ is not a root of unity then
the center of $\cS(F_{0,4};R,A)$ is the subalgebra generated by boundary
components.
\end{theorem}

\begin{proof}
We use the the complexity $a+b+c$ for $(a,b,c)$.  Suppose $\alpha$ is
a central element.  Write it as $\sum r_i (a_i,b_i,c_i) + \beta$,
where each $(a_i,b_i,c_i)$ is distinct but all have the same
complexity, and $\beta$ consists of lower complexity terms.  Lemma 4
implies $0=[\alpha,(2,0,0)] = \sum r_i(A^{c_i-b_i} - A^{b_i-c_i})
(a_i+2,b_i,c_i) + \mathcal{O}$, so $b_i = c_i$.  Similar computations with
$(0,2,0)$ and $(0,0,2)$ show that each $(a_i,b_i,c_i)$ is a power of
the outer boundary curve.  That means $\beta$ is central, so we can
repeat the same argument for it.  Continuing in this fashion, we see
that every term of $\alpha$ is a scalar times a power of the outer
boundary component.
\end{proof}

\begin{theorem}
If $R$ has no zero divisors then neither does $\cS(F_{0,4};R,A)$. 
\end{theorem}

\begin{proof}
We use the complexity $a^2+b^2+c^2$ for $(a,b,c)$.  Choose $\alpha =
\sum r_i(a_i,b_i,c_i)$ with $r_i \neq 0$ and each $(a_i,b_i,c_i)$
distinct.  Similarly choose $\beta = \sum s_j(a'_j,b'_j,c'_j)$.
Reordering if necessary, we may assume $(a_1+a'_1,b_1+b'_1,c_1+c'_1)$
is a maximal complexity link in the resolution of $\alpha \beta$.  By
Lemma 4, it occurs in the resolution of
$r_1s_1(a_1,b_1,c_1)(a'_1,b'_1,c'_1)$ with non-zero coefficient.  By
Lemma 3 it does not occur in any other
$r_is_j(a_i,b_i,c_i)(a'_j,b'_j,c'_j)$.  Hence, $\alpha \beta \neq 0$.
\end{proof}

The central element theorems can fail at roots of unity.  When $A^2
=1$ any $\cS(F;R,A)$ is commutative.  When $A^4=1$, $\cS(F_{0,4};R,A)$
is still commutative.  The toral examples are not, but the square of
every curve is central.

\end{document}